\documentclass[sigconf]{acmart}

\usepackage{booktabs} 

\usepackage{amsmath}
\usepackage[caption=false,font=footnotesize]{subfig}
\usepackage{url}

\usepackage{graphicx}
\usepackage{amsfonts}
\usepackage{amssymb}
\usepackage{amsthm}
\usepackage{hyperref}
\usepackage{color}
\usepackage{colortbl}
\usepackage{mydef}

\setcopyright{rightsretained}

\acmDOI{http://dx.doi.org/10.1145/3067695.3076069}

\acmISBN{978-1-4503-4939-0/17/07}

\acmConference{GECCO '17 Companion}{July 15-19, 2017}{Berlin, Germany}
\acmYear{2017}
\copyrightyear{2017}

\acmPrice{15.00}

\begin{document}
\title{Simple Problems}
\subtitle{The Simplicial Gluing Structure of Pareto Sets and Pareto Fronts}

\author{Naoki Hamada}
\affiliation{%
  \institution{Fujitsu Laboratories Ltd.}
  \streetaddress{4-1-1 Kamikodanaka, Nakahara-ku}
  \city{Kawasaki} 
  \country{Japan}}
\email{hamada-naoki@jp.fujitsu.com}

\renewcommand{\shortauthors}{N. Hamada}

\begin{abstract}
Quite a few studies on real-world applications of multi-objective optimization reported that their Pareto sets and Pareto fronts form a topological simplex.
Such a class of problems was recently named the \emph{simple problems}, and their Pareto set and Pareto front were observed to have a gluing structure similar to the faces of a simplex.
This paper gives a theoretical justification for that observation by proving the gluing structure of the Pareto sets/fronts of subproblems of a simple problem.
The simplicity of standard benchmark problems is studied.
\end{abstract}

%
%

\ccsdesc[500]{Applied computing~Multi-criterion optimization and decision-making}
\ccsdesc[500]{Mathematics of computing~Evolutionary algorithms}
\ccsdesc[500]{Mathematics of computing~Nonconvex optimization}
\ccsdesc[500]{Mathematics of computing~Geometric topology}


\keywords{multi-objective optimization, continuous optimization, problem class, stratification}

\maketitle

\section{Introduction} \label{sec:introduction}
\subsection{Motivation}
The success of evolutionary multi-objective optimization (EMO) is widely spreading over various academic and industrial fields.
Recent numerical studies showed that decomposition-based EMO algorithms such as MOEA/D~\cite{Zhang07}, NSGA-III~\cite{Deb14}, and AWA~\cite{Hamada10,Hamada11a,Hamada11b,Shioda15,Shioda15b} have an ability to approximate the entire Pareto set and Pareto front of many-objective problems.

In contrast to their abundance of experimental successes, the theory shedding light on why they work is still under developing.
Especially, the problem class in which decomposition-based EMO algorithms can cover the entire Pareto set/front has not been understood.
This paper discusses some problem class in which solutions are well-behaved for scalarization.

\subsection{Related Studies in Other Fields}
Theories concerning the easiness of covering solutions have been developed in several optimization-related fields as well as the EMO community.
Most of them are studies on topological properties of solution sets.

\subsubsection{Contractibility}
The earliest work can be found in 1951 by Koopmans~\cite{Koopmans51} Assertion 4.14 in which he applied the linear programming to economics and pointed out some conditions making the Pareto front contractible.
Peleg~\cite{Peleg72} generalized this result and showed that the Pareto front is contractible if the feasible objective region is a convex set.
Afterward, the study spread to operations research, and the closedness, the (arcwise) connectedness, and the contractibility of the Pareto set/front under general settings were studied on linear programming in 1970's and on (quasi) convex programming in 1980's.
These results are collected in Luc~\cite{Luc89} Section 6.
Recently, similar results were obtained under more general problem classes such as the lexicographic quasiconvexity~\cite{Popovici06} and the arcwise cone-quasiconvexity~\cite{LaTorre10}.
The (arcwise) connectedness of the Pareto set is a necessary condition that the homotopy method covers the solutions.

\subsubsection{Decomposition}
The decomposition approach considers not only a given problem but also its subproblems each optimizing a subset of objective functions and studies the relation among their solutions.
Lowe et~al.~\cite{Lowe84} in 1984 showed that the weak Pareto set of a convex programming is the union of the Pareto sets of subproblems.
Malivert et~al.~\cite{Malivert94} extended the result to explicitly quasi-convex upper semicontinuous functions.
Popovici~\cite{Popovici05} named this property the \emph{Pareto reducible} and gave a sufficient condition independent of convexity.
Ward~\cite{Ward89} showed that the strictly Pareto solutions to a convex programming problem are \emph{completely surrounded} by the Pareto solutions of subproblems.
Recent studies~\cite{Popovici06,Popovici07,LaTorre10} revealed that the Pareto reducibility of the lexicographic quasiconvex programming problem is closely related to the contractibility, through the \emph{simply shadiness}~\cite{Benoist03} of the Pareto front.

\subsubsection{Stratification}
From pure mathematics, the singularity theory of differentiable maps gives a decomposition of solutions.
In 1973, Smale~\cite{Smale73} applied this theory to an economic problem and stated that the Pareto set of a pure exchange economy of $m$ agents is homeomorphic to an $(m-1)$-simplex, provided the quasiconvexity and monotonicity of the agents' utility functions.
Lovison et~al.~\cite{Lovison14} pointed out that each face of this simplex corresponds to the Pareto set of each subproblem optimizing a subset of objective functions.
In Smale~\cite{Smale73} and its sequels~\cite{Smale74a,Smale74b,Smale74c,Smale74d,Smale76}, he discussed the \emph{stratification} of Pareto critical points of generic maps with the transversality and rank assumption of derivatives.
de Melo~\cite{deMelo76} showed that the $C^\infty$-maps whose Pareto critical points admit a stratification are generic, i.e., they form a dense subset of the space of $C^\infty$-maps under the Whitney topology.
Recently, Lovison et~al.~\cite{Lovison14} collected related works to this approach and further developed de Melo's result, and showed that local Pareto sets of sufficiently proper maps admit a \emph{Whitney stratification}.

\subsection{Our Approach}
These attempts are going on in two courses: the linear/convex analysis originated by Koopmans and the global analysis by Smale.
The former seems much restrictive for the global optimization nature of EMO algorithms.
The latter approach is sufficiently general but currently hard to compute.
We need a handy theory for understanding the behavior of EMO algorithms.
Recently, Hamada et~al.~\cite{Hamada11a} defined a class of problems called the \emph{simple problem}.
They pointed out, without rigorous proofs, that the Pareto set and Pareto front of a simple problem are both homeomorphic to a simplex and the faces of the simplex correspond to the Pareto sets and their images of the subproblems.
They also discussed that this property is closely related to scalarization.
This paper gives rigorous proofs for their arguments.

\subsection{Contribution}
In this paper, we give a proof that the boundary of the Pareto set (resp.~Pareto front) of a simple problem is the union of the interior of the Pareto sets (resp.~their images) of subproblems. This property enables us to numerically compute a stratification of the Pareto set/front.

Additionally, we investigate the simplicity of benchmark problems widely-used in the EMO community: all problems in ZDT suite~\cite{Zitzler00} and DTLZ suite~\cite{Deb05b} are non-simple, five of problems in WFG suite~\cite{Huband06} can be simple under a very restrictive situation, and MED problem~\cite{Hamada11b} is always simple.

\subsection{Contents}
The rest of the paper is organized as follows.
Section \ref{sec:preliminaries} prepares basic notions and notations used in subsequent sections.
Section \ref{sec:main results} gives some properties of solutions to simple problems and their relation to scalarization.
Section \ref{sec:benchmark} discusses the simplicity of existing benchmark problems.
Section \ref{sec:conclusions} gives conclusions and remarks for future work.

\section{Preliminaries} \label{sec:preliminaries}
This paper considers the following \emph{$n$-variable $m$-objective problem}:
\begin{equation} \label{eqn:MCOP}
 \minimize_{x \in X \subseteq \R^n} \f(x) = (f_1(x),\ldots,f_m(x))
\end{equation}
where we call $\f: \R^n \to \R^m$ the \emph{evaluation map}, $f_i:\R^n \to \R \ (i = 1, \ldots, m)$ the \emph{$i$-th objective function}, $\R^n$ the \emph{variable space}, $X \subseteq \R^n$ the \emph{feasible region}, $x \in \R^n$ a \emph{solution}, $\R^m$ the \emph{objective space}, and $\f(x) \in \R^m$ an \emph{evaluation value}.

We will make various problems by removing some of objective functions from \eqref{eqn:MCOP} and discuss a gluing structure of their solutions.
To write such arguments clearly, we define the problem as a finite set of objective functions $F = \set{f_1, \dots, f_m}$ and regard \eqref{eqn:MCOP} as a notation for it.
We abuse $\f$ to denote the problem $F$, the equation \eqref{eqn:MCOP}, and the evaluation map $\f$, depending on the context.
We call the empty set $\emptyset$ a \emph{$0$-objective problem} and define its evaluation map as the empty map $\zeroproblem: \emptyset \to \R^0$ to a one-point set $\R^0 = \set{0}$.

For problems $\f, \g$ such that $\g \subproblemeq \f$ as sets, we say that $\g$ is a \emph{subproblem} of $\f$ and $\f$ is a \emph{superproblem} of $\g$.
We call the set of all subproblems of a problem $\f$ the \emph{decomposition} of $\f$ and denote by $\decomp{\f} = \set{\g \mid \g \subproblemeq \f}$.

If solutions $x, y \in X$ to a problem $\f$ satisfy the conditions
\[
\forall f_i \inproblem \f,\ f_i(x) \le f_i(y) \ \text{ and } \ \exists f_j \inproblem \f,\ f_j(x) < f_j(y),
\]
then we say that $x$ \emph{$\f$-dominates} $y$ and denote it by $x \dominates_{\f} y$.
We denote $\f(x) = \f(y)$ by $x =_{\f} y$ and $(x \dominates_{\f} y) \lor (x =_{\f} y)$ by $x \dominateseq_{\f} y$.

If a solution $x^* \in X$ to a problem $\f$ satisfies
\[
\forall x \in X,\ x \not \dominates_{\f} x^*,
\]
then $x^*$ is called a \emph{Pareto solution} to $\f$.

The set of all Pareto solutions of a problem $\f = \set{f_1, \dots, f_m}$ is called the \emph{Pareto set}, denoted by $X^*(\f)$ or $X^*(f_1,\ldots,f_m)$.
The image of $X^*(\f)$ by a map $\g$ is denoted by $\g X^*(\f)$.
Especially, the image $\f X^*(\f)$ is called the \emph{Pareto front}.
Through the paper, we abbreviate the composition of maps $\psi \circ \phi$ to $\psi \phi$.
The above notation, $\g X^*$, can be considered as $\g \circ X^*$ by regarding $X^*$ as a map $X^*: 2^F \to 2^X$.

\section{Simple Problem} \label{sec:main results}
This section introduces the definition of simple problem and shows its solution structure.
Section~\ref{sec:simple problem} presents the definition of simple problem.
Section~\ref{sec:solution properties} shows some inclusion properties of solutions among subproblems.
Section~\ref{sec:gluing properties} shows that those solutions have the gluing structure of a topological simplex.
Section~\ref{sec:AWA} points out that this gluing structure enables decomposition-based EMO algorithms to cover the Pareto set and Pareto front.

\subsection{Definition} \label{sec:simple problem}
First, we present the definition of simple problem and its graphical intuition.
\begin{figure*}[t]
\begin{tabular}{ccc}
\mpg[0.3]{\tiny
\begin{eqnarray*}
\minimize_{x_1,x_2 \in [-2,2]}\\
f_1(x) =& \hspace{-2em} x_1^2 &+\ 3(x_2 - 1)^2\\
f_2(x) =& 2(x_1 - 1)^2             &+\ \quad x_2^2\\
f_3(x) =& 3(x_1 + 1)^2             &+\ 2(x_2 + 1)^2
\end{eqnarray*}}
\ig[0.3]{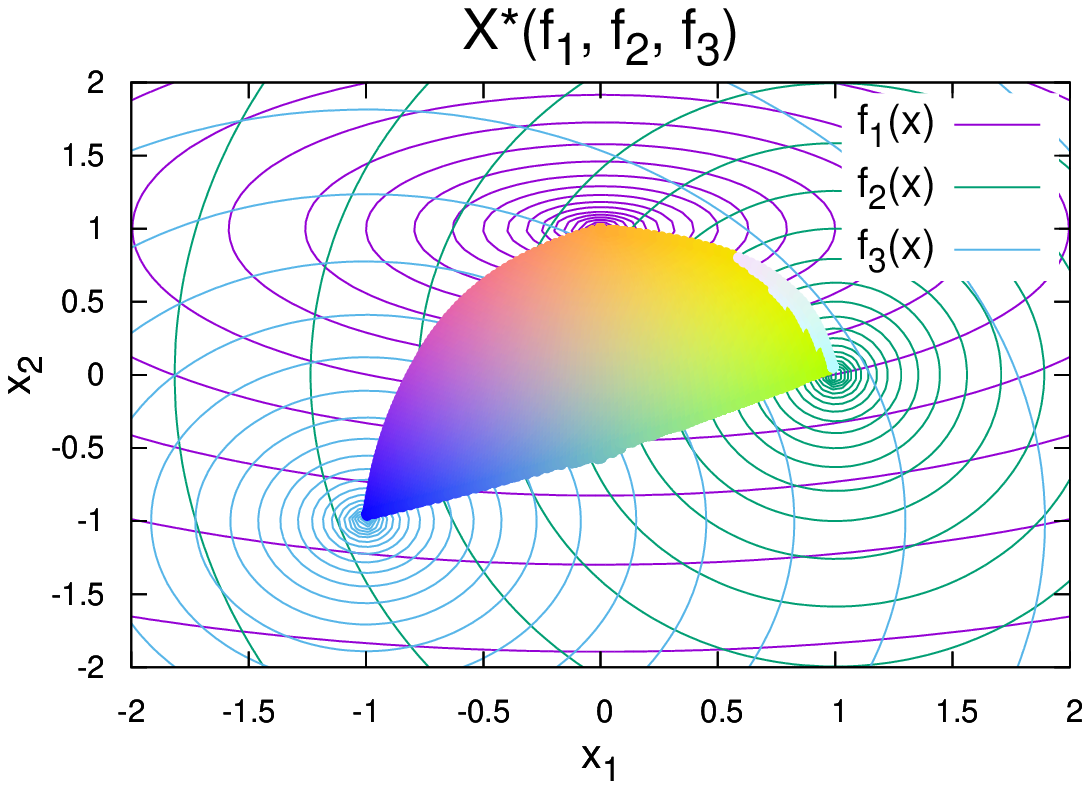}
\mpg[0.3]{\tiny
\begin{eqnarray*}
\mathrm{RGB}(x) &=& (y_1(x), y_2(x), y_3(x))\\
y_i(x) &=& \begin{cases}
 \paren{\ol{f_i} - f_i(x)} / \paren{\ol{f_i} - \ul{f_i}} & (f_i \in \g)\\
 0 & (f_i \not \in \g)
\end{cases}\\
\ol{f_i} &=& \max \set{f_i(x) \mid x \in X^*(\g)}+\mbox{1e-16}\\
\ul{f_i} &=& \min \set{f_i(x) \mid x \in X^*(\g)}
\end{eqnarray*}}\\
\ig[0.3]{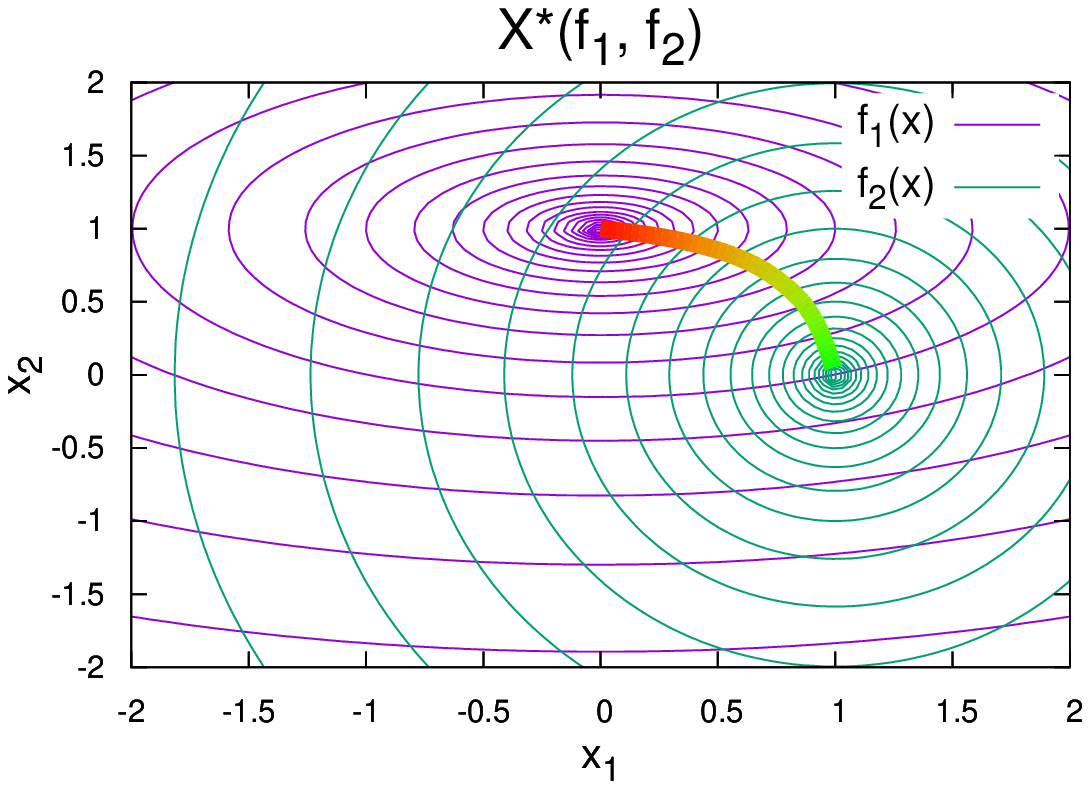}
\ig[0.3]{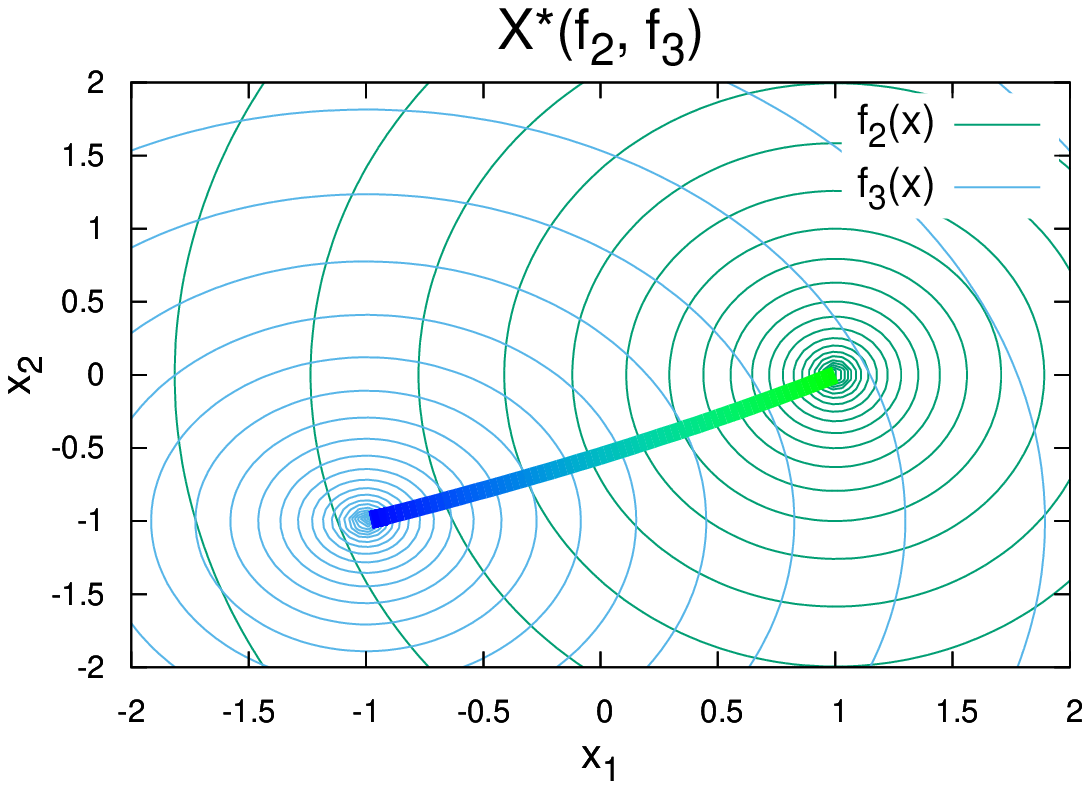}
\ig[0.3]{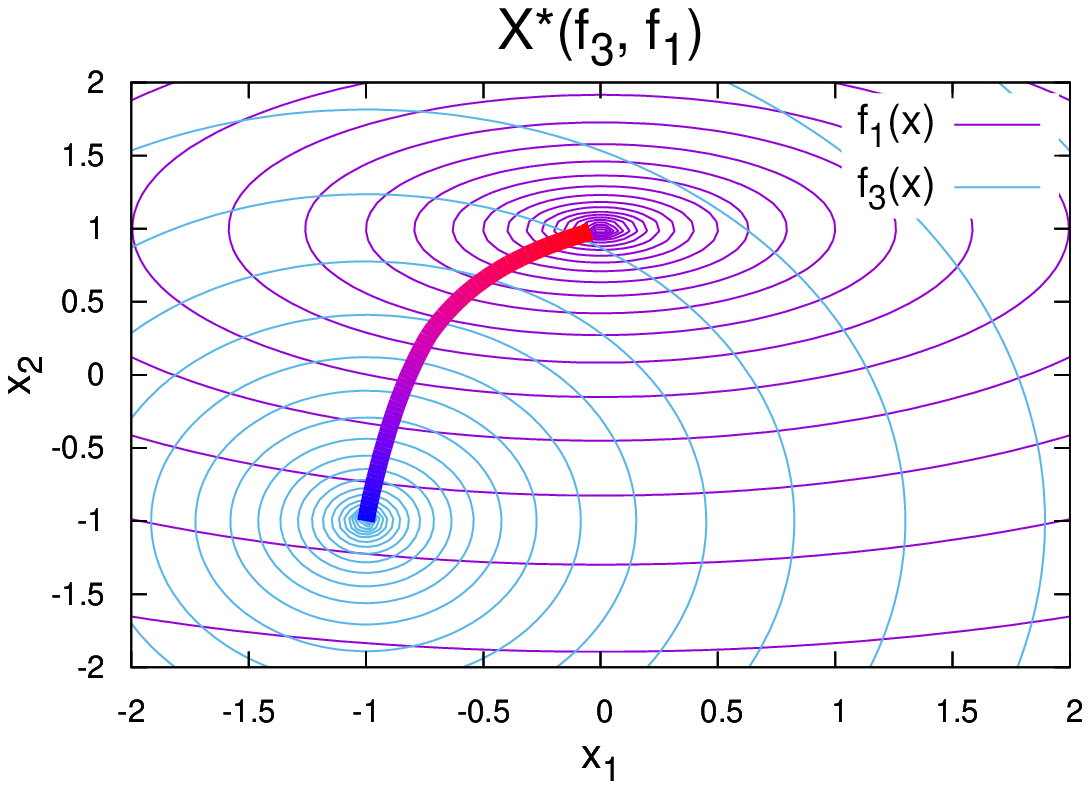}\\
\ig[0.3]{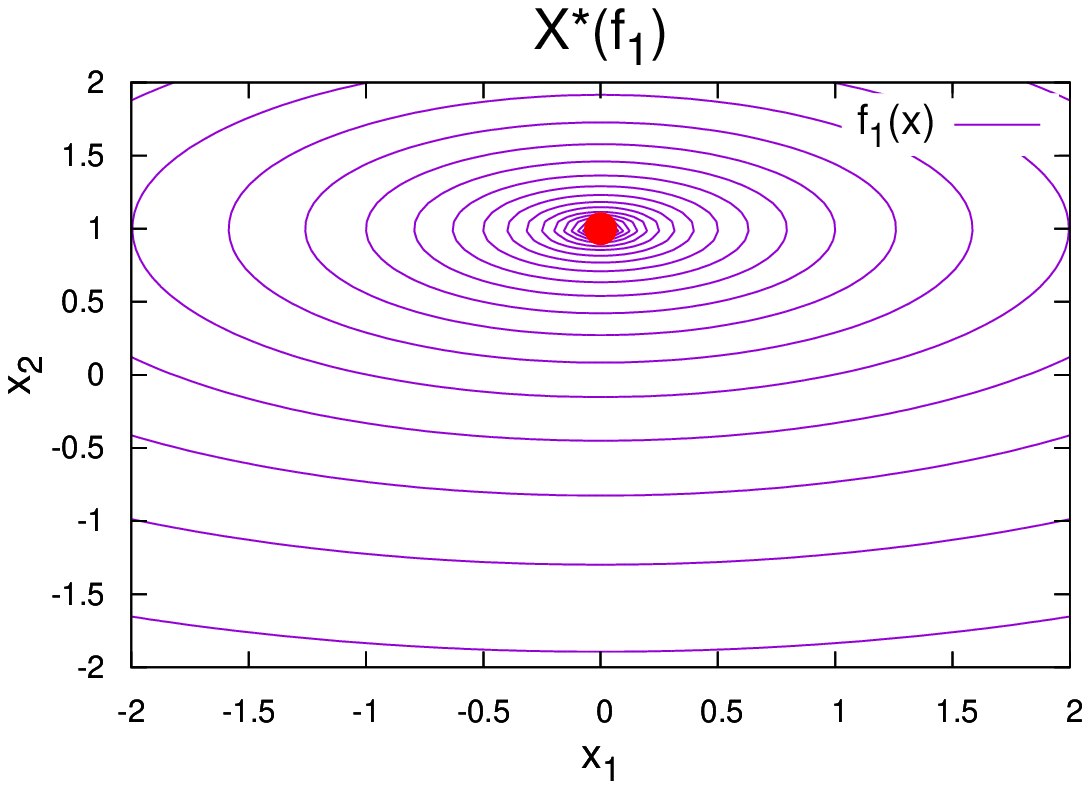}
\ig[0.3]{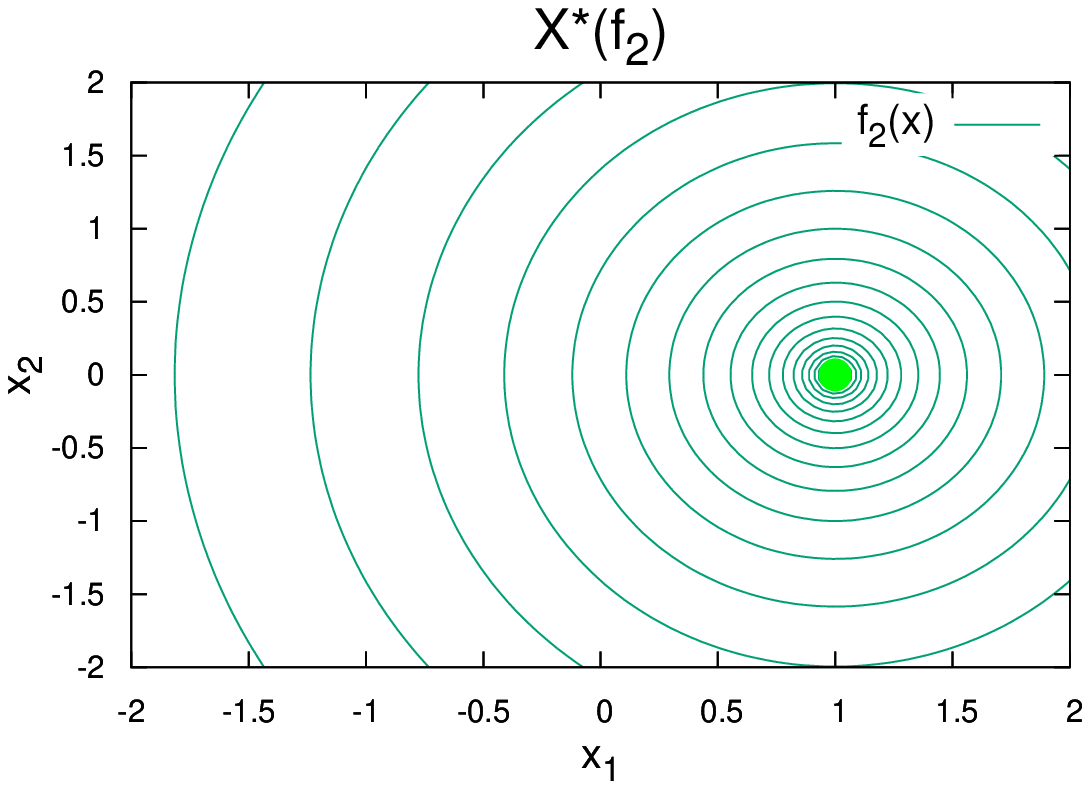}
\ig[0.3]{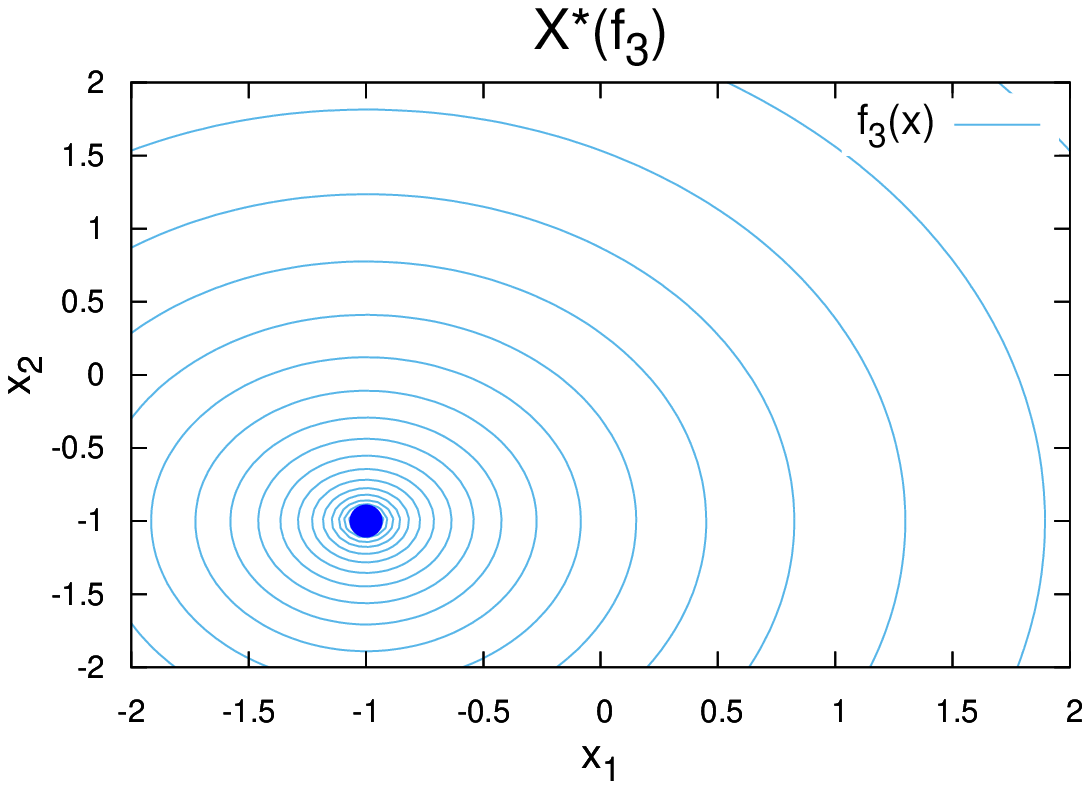}\\
\mpg{\caption{A simple 3-objective problem $\f = \set{f_1,f_2,f_3}$ and its subproblem $\g$'s Pareto set $X^*(\g)$. The Pareto sets are colored by converting the $f_1,f_2,f_3$-coordinates to RGB using the equations in the top-right cell.}
\label{fig:simple}}
\end{tabular}
\end{figure*}

\begin{defn}[simple problem~\cite{Hamada11a}] \label{def:simplicity}
A problem $\f$ is \emph{simple} or  has \emph{simplicity} if every subproblem $\g \in 2^{\f}$ satisfies the following conditions:
\begin{description}
 \item[(S1)] $X^*(\g) \homeo \Delta^{k-1}$ if $\g$ is $k$-objective,
 \item[(S2)] $\g|_{X^*(\g)}: X^*(\g) \to \R^k$ is an embedding.
\end{description}
Here, $\Delta^{k-1} = \set{(x_1, \dots, x_k) \in [0,1]^k \mid \sum x_i =1}$ is the \emph{standard $(k-1)$-simplex}.
By $A \homeo B$, we denote that topological spaces $A$ and $B$ are \emph{homeomorphic}, which is defined as there are continuous maps $\psi: A \to B$ and $\phi: B \to A$ such that $\phi \psi = \mathrm{id}_A$ and $\psi \phi = \mathrm{id}_B$.
Such maps, $\psi$ and $\phi$, are called \emph{homeomorphisms}.
The topology of $X^*(\g)$ is induced from the variable space $\R^n$.
That is, any open set $U$ in $X^*(\g)$ can be written as $U = V \cap X^*(\g)$ with some open set $V$ in $\R^n$ under the Euclidean topology.
Similarly, all spaces discussed in this paper are implicitly equipped with induced topologies from either the variable space $\R^n$ or the objective space $\R^m$.
By $\phi|_{A'}: A' \to B$, we denote the \emph{restriction} of a map $\phi: A \to B$ to a set $A' \subseteq A$, which is the composite $\phi \iota: A' \to B$ with the inclusion map $\iota: A' \embedding A, a \mapsto a$.
The \emph{embedding} of $A$ to $B$ is the composite $\iota \phi: A \to B$ of a homeomorphism $\phi: A \to B'$ and the inclusion map $\iota: B' \embedding B$.
\end{defn}
Let us cultivate the intuitive understanding of this definition with Fig.~\ref{fig:simple}.
When considering a simple problem, we also deal with its subproblems optimizing all subsets of given objective functions.
The first subproblem is the original problem itself optimizing all the objectives.
The example in the figure is $3$-objective, and the condition (S1) imposes that its Pareto set, $X^*(f_1,f_2,f_3)$, is homeomorphic to $\Delta^2$, a surface created by bending and stretching a triangle (without cutting and connecting).
The condition (S2) guarantees that the restricted evaluation map $\f: X^*(f_1, f_2, f_3) \to \f X^*(f_1, f_2, f_3)$ and its inverse map $\f^{-1}: \f X^*(f_1, f_2, f_3) \to X^*(f_1, f_2, f_3)$ are bijective, continuous, and thus homeomorphisms.
This implies that the Pareto front, $\f X^*(f_1, f_2, f_3)$, is also homeomorphic to $\Delta^2$, and every point on $\f X^*(f_1, f_2, f_3)$ continuously corresponds to a unique solution on $X^*(f_1, f_2, f_3)$ and vice versa.
Next, we remove one of objectives, resulting in three 2-objective subproblems.
Their Pareto sets, $X^*(f_1,f_2)$, $X^*(f_2,f_3)$, $X^*(f_3,f_1)$, are homeomorphic to $\Delta^1$ by (S1), a curve without loops.
By (S2), their Pareto fronts are also a curve in which each point continuously corresponds to a unique Pareto solution.
Again removing an objective, we get three 1-objective subproblems.
Their Pareto sets, $X^*(f_1)$, $X^*(f_2)$, $X^*(f_3)$, are homeomorphic to $\Delta^0$, a point!
A map on a one-point set always satisfies (S2); there is no special implication.
Finally, there is a 0-objective problem (not shown in Fig.~\ref{fig:simple}) that corresponds to the case of no objective function.
This is just required for formality.

Note that the conditions (S1) and (S2) only impose the problem structure within the Pareto set and admit an arbitrary structure out of the Pareto set.
This is contrastive to conventional problem classes such as linear/convex/polynomial programming problems which regulate their structures over the entire domain.
The simplicity is independent of those problem classes.
In fact, every class of linear/convex/polynomial programming problems contains both of simple problems and non-simple problems!
As a result, the simplicity characterizes a new aspect of ``easiness to solve''.

\subsection{Inclusion Properties} \label{sec:solution properties}
This section shows inclusion relations of the Pareto sets and its images for subproblems of a simple problem.
To isolate the consequence of assuming the simplicity from general properties of Pareto-optimality, we begin without the assumption.
\begin{prp} \label{prp:subproblem}
For any problem $\f$ (possibly non-simple) and any subproblem $\g \subproblemeq \f$, the following relations hold:
\[
X^*(\g) \subseteq \wX(\g) \subseteq \wX(\f)
\]
where $\wX(\g)$ is the \emph{weak Pareto set} of $\g$, that is, the set of points each $x \in X$ satisfying
\begin{equation} \label{eqn:weak Pareto}
\forall y \in X, \exists f_i \in \g: f_i(x) \le f_i(y).
\end{equation}
\end{prp}
\begin{prf}
The first inclusion, $X^*(\g) \subseteq \wX(\g)$, is a well-known fact (see, for example, Miettinen~\cite{Miettinen99} Section~2.5).
The second relation, $\wX(\g) \subseteq \wX(\f)$, directly follows from \eqref{eqn:weak Pareto}.
\end{prf}
However, it does not hold that $X^*(\g) \subseteq X^*(\f)$ in general.
\begin{eg} \label{eg:non pareto}
Consider a 1-variable 2-objective problem
\[
\minimize_{x \in [0,1]} f_1(x)=0, f_2(x)=x.
\]
Clearly, the Pareto sets are $X^*(f_1) = [0,1]$ and $X^*(f_1, f_2) = \set{0}$, which implies $X^*(f_1) \not \subseteq X^*(f_1,f_2)$.
\end{eg}

In contrast, simple problems do not have such ill-behaved solutions.
\begin{prp} \label{prp:no WP}
If a problem $\f$ is simple, then
\[
\wX(\f) = X^*(\f).
\]
\end{prp}
\begin{prf}
From Proposition~\ref{prp:subproblem}, we have $\wX(\f) \supseteq X^*(\f)$.
We will prove $\wX(\f) \subseteq X^*(\f)$ by contradiction.
Suppose that a point $x \in \wX(\f) \setminus X^*(\f)$ exists.
Then, since $x$ is weakly Pareto-optimal, the condition
\begin{equation} \label{eqn:WP}
\forall y \in X, \exists f_i \in \f: f_i(x) \le f_i(y)
\end{equation}
holds.
On the other hand, since $x$ is not Pareto-optimal, the condition
\begin{equation} \label{eqn:NP}
\exists y \in X \setminus \set{x}, \forall f_i \in \f: f_i(y) \le f_i(x)
\end{equation}
holds.
A point $y$ in \eqref{eqn:NP} should satisfy \eqref{eqn:WP}, which means
\[
\exists f_i \in \f: f_i(y) = f_i(x).
\]
By \eqref{eqn:WP}, we have $x \in X^*(f_i)$ and thus $y \in X^*(f_i)$.
This contradicts (S2); more specifically, $f_i$ cannot be an injection on $X^*(f_i)$.
\end{prf}
Therefore, the simplicity ensures the inclusion relationship of Pareto sets.
\begin{prp} \label{cor:inclusion}
For a simple problem $\f$ and any subproblem $\g \subproblemeq \f$, it holds that
\[
X^*(\g) \subseteq X^*(\f).
\]
\end{prp}
\begin{prf}
Combine Proposition~\ref{prp:subproblem} with Proposition~\ref{prp:no WP}.
\end{prf}

Using this fact, we can see the topology of the image of the Pareto set.
Although, the condition (S1) itself addresses only the topology of the Pareto set, combined with (S2), we show that the Pareto front and the image in the superproblem have the same topology.
\begin{prp} \label{prp:simple homeo}
For a simple problem $\f$ and any subproblem $\g \subseteq \f$, if $\g$ is $k$-objective, then
\[
X^*(\g) \homeo \g X^*(\g) \homeo \f X^*(\g) \homeo \Delta^{k-1}.
\]
\end{prp}
\begin{prf}
By Definition~\ref{def:simplicity}, $X^*(\g) \homeo \g X^*(\g) \homeo \Delta^{k-1}$ is evident.
We will show the nontrivial part: $X^*(\g) \homeo \f X^*(\g)$.
It suffices to show that the restriction $\f|_{X^*(\g)}$ is an embedding.
Remember that the restriction $\f|_{X^*(\f)}$ is an embedding by definition and $X^*(\g) \subseteq X^*(\f)$ in Proposition~\ref{cor:inclusion}.
In general, any restriction of an embedding is again an embedding.
\end{prf}

In the definition of simplicity, the conditions (S1) and (S2) are imposed on all subproblems as well as a given problem, which means that the subproblems inherit the simplicity from the superproblem.
\begin{prp} \label{lem:simple}
If a problem $\f$ is simple, then any subproblem $\g \subproblemeq \f$ is simple.
\end{prp}
\begin{prf}
If $\g \subproblemeq \f$, then $2^{\g} \subseteq 2^{\f}$.
Thus, if all the problems in $2^{\f}$ satisfy (S1) and (S2), then also do in $2^{\g}$, which implies that if $\f$ is simple, then $\g$ is.
\end{prf}
Therefore, propositions that hold for a simple problem also hold for its subproblems.
For example, the actual assertion of Proposition~\ref{prp:no WP} is that any subproblem of $\f$, as well as $\f$ itself, does not have a weakly Pareto-optimal point which is not Pareto-optimal.
The interpretation of Proposition~\ref{prp:simple homeo} is bit more complicated: given a simple problem $\h$, the assertion holds for any pair of problems $\f, \g$ such that $\g \subseteq \f \subseteq \h$.
Henceforth, we will not repeat this property, but it is always valid when propositions involve simple problems.

Similarly to that the empty set is a subset of every set, a 0-objective problem is a subproblem of every problem.
Therefore, if 0-objective problems do not exist or are not simple, then there cannot exist any simple problem.
Let us check them.
\begin{prp} \label{lem:simple0}
There exists a $0$-objective problem $\zeroproblem$; it is unique and simple.
\end{prp}
\begin{prf}
The existence and the uniqueness follow from those of the empty set and empty map.
Let us check the simplicity.
Since the only subproblem of $\zeroproblem$ is $\zeroproblem$ itself, it suffices to show that $\zeroproblem$ satisfies (S1) and (S2).

(S1):
Since $\zeroproblem$ is the 0-objective problem, it suffices to show $X^*(\zeroproblem) \homeo \Delta^{-1}$.
It holds indeed as $\Delta^{-1} = \emptyset$ and $X^*(\zeroproblem) = \emptyset$.

(S2):
By $X^*(\zeroproblem) = \emptyset$, the restricted evaluation map, $\zeroproblem|_{X^*(\zeroproblem)}: X^*(\zeroproblem) \to \R^0$, is the evaluation map $\zeroproblem: \emptyset \to \R^0$ itself.
Since $\zeroproblem: \emptyset \to \R^0$ can be decomposed into a homeomorphism $\mathrm{id}: \emptyset \to \emptyset$ and an inclusion map $\iota: \emptyset \embedding \R^0$ as $\zeroproblem = \iota \mathrm{id}$, the restriction $\zeroproblem|_{X^*(\zeroproblem)}: X^*(\zeroproblem) \to \R^0$ is an embedding.
\end{prf}

Now we have confirmed that for each subproblem, the Pareto set and its image are well-behaved.
The next section investigates that those solutions are nicely glued together.

\subsection{Gluing Properties} \label{sec:gluing properties}
\begin{figure}[t]
\includegraphics[width=\hsize]{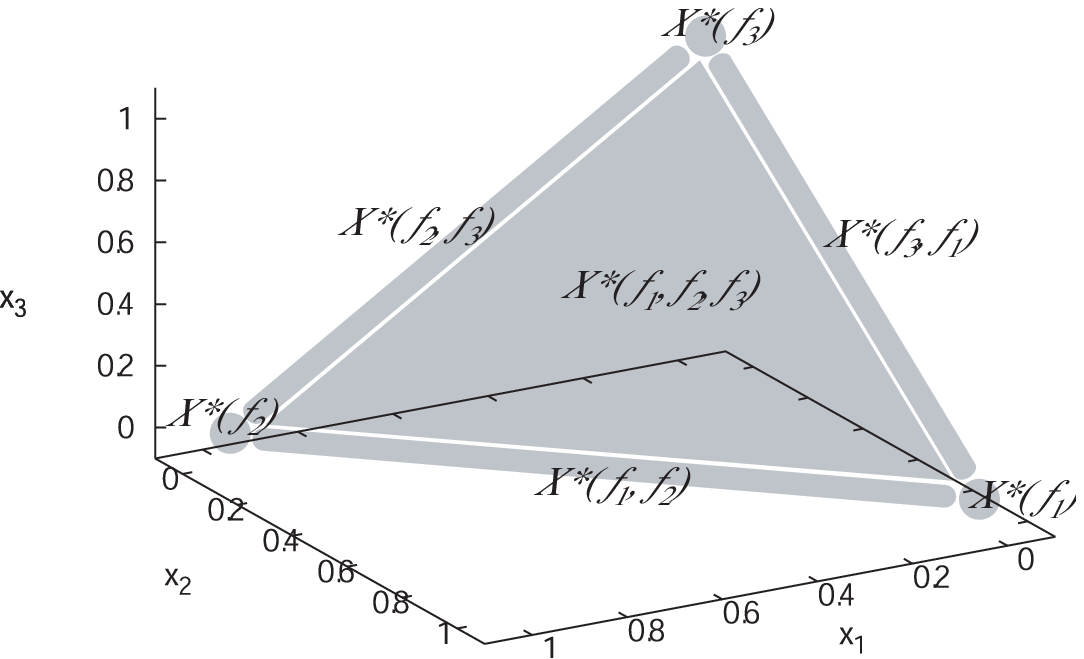}
\caption{A simple 3-objective problem $\f = \set{f_1,f_2,f_3}$ and its gluing structure of the Pareto sets $X^*(\g)$ of subproblems $\g \subproblemeq \f$. The boundary of the Pareto set $X^*(\g)$ of each subproblem $\g \subproblemeq \f$ consists of the Pareto sets $X^*(\h)$ of all subproblems $\h \subproblem \g$.}
\label{fig:simple problem X}
\includegraphics[width=\hsize]{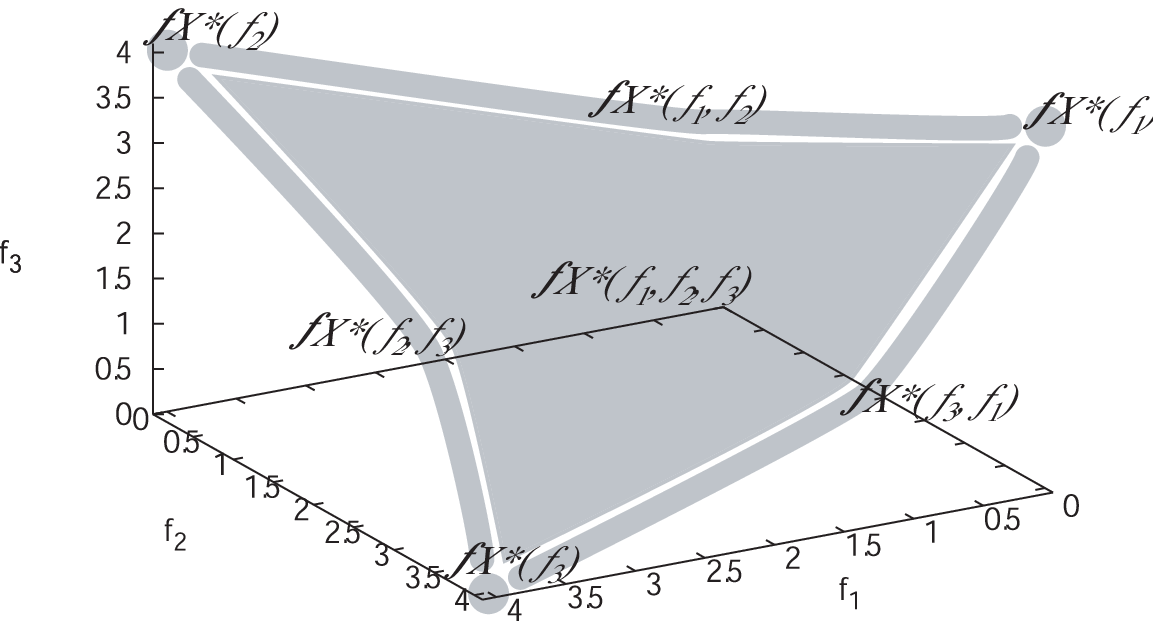}
\caption{A simple 3-objective problem $\f = \set{f_1,f_2,f_3}$ and its gluing structure of the Pareto set images $\f X^*(\g)$ for all subproblems $\g$. Although the shape is different from Fig.~\ref{fig:simple problem X}, the topology is the same.}
\label{fig:simple problem F}
\end{figure}

The goal of this section is to give a proof that the solutions to a simple problem have a special gluing structure as shown in Figs.~\ref{fig:simple problem X} and \ref{fig:simple problem F}.
This structure is an analogy of the faces of a simplex.
A $2$-simplex $[v_1, v_2, v_3]$ is a triangle spanned by vertices $v_1, v_2, v_3$ and its boundary $\boundary[v_1, v_2, v_3]$ is the union of three edges, or $1$-simplices, $[v_1, v_2], [v_2, v_3], [v_3, v_1]$.
Each edge $[v_i, v_j]$ has the boundary consisting of two points, or $0$-simplices, $[v_i]$ and $[v_j]$.
The boundary of each vertex $[v_k]$ is the empty set, or the $(-1)$-simplex.
We can expand these relations using $A = \Int A \sqcup \boundary A$ and see that the boundary of a simplex is expressed as the disjoint union of the open faces:
\begin{multline*}
\boundary [v_1, v_2, v_3] =\\
 \Int [v_1, v_2] \sqcup \Int [v_2, v_3] \sqcup \Int [v_3, v_1]\\
 \sqcup \Int [v_1] \sqcup \Int [v_2] \sqcup \Int [v_3].
\end{multline*}

Generally, Pareto sets and their images may have a more complex topological structure.
However, for a simple problem, we have seen in Proposition \ref{prp:simple homeo} that those form topological manifolds with boundary (hereafter, simply call \emph{manifolds}) homeomorphic to a simplex.
In a $k$-manifold $M$, a point having an open neighborhood homeomorphic to $\R^k$ is called an \emph{interior point}, and the set of all interior points is called the \emph{interior} of $M$, denoted by $\Int M$.
The other points are \emph{boundary points} having an open neighborhood homeomorphic to $[0, \infty) \times \R^{k-1}$, and the set of all boundary points is called the \emph{boundary} of $M$, denoted by $\boundary M$.
For a simple problem, a similar relation holds among Pareto sets as shown in Fig.~\ref{fig:simple problem X}:
\begin{multline} \label{eqn:gluing X}
\boundary X^*(f_1, f_2, f_3) =\\
 \Int X^*(f_1, f_2) \sqcup \Int X^*(f_2, f_3) \sqcup \Int X^*(f_3, f_1)\\
 \sqcup \Int X^*(f_1) \sqcup \Int X^*(f_2) \sqcup \Int X^*(f_3).
\end{multline}
The same relation holds for the images as shown in Fig.~\ref{fig:simple problem F}:
\begin{multline} \label{eqn:gluing fX}
\boundary \f X^*(f_1, f_2, f_3) =\\
 \Int \f X^*(f_1, f_2) \sqcup \Int \f X^*(f_2, f_3) \sqcup \Int \f X^*(f_3, f_1)\\
 \sqcup \Int \f X^*(f_1) \sqcup \Int \f X^*(f_2) \sqcup \Int \f X^*(f_3).
\end{multline}
It is known that such a gluing structure of solutions commonly appears in facility location problems, studied for a long time (see for example Rodr\'{i}guez-Ch\'{i}a et~al.~\cite{Rodriguez-Chia02} and the references therein).
It is also seen in other applications and exploited as a heuristic for practitioners \cite{Everson13}.

Now we start to show \eqref{eqn:gluing X} and \eqref{eqn:gluing fX} for an arbitrary number of objectives.
First of all, let us see some basic properties of the Pareto front that hold for any class of problems.
\begin{lem} \label{lem:projection}
For an $m$-objective (possibly non-simple) problem $\f$ whose Pareto front $\f X^*(\f)$ forms an $(m-1)$-manifold,
the projection
\[
\pi_{-i}: \left\{
\begin{array}{ccc}
\R^m & \to    & \paren{\R^{i-1} \times \set{ 0 } \times \R^{m-i}} \simeq \R^{m-1},\\
(y_1, \dots, y_m)     &\mapsto & (y_1, \dots, y_{i-1}, 0, y_{i+1}, \dots, y_m)
\end{array}
\right.
\]
restricted to $\Int \f X^*(\f)$ is an embedding.
\end{lem}
\begin{prf}
Generally, any projection and its restriction to any open set are continuous and open, and any injective continuous open map is an embedding.
Thus, it suffices to show that $\pi_{-i}$ is injective when restricted to $\Int \f X^*(\f)$.
If $\pi_{-i}|_{\Int \f X^*(\f)}$ is not injective, then $\Int \f X^*(\f)$ contains two points having the same coordinates except for the $i$-th value.
This means that one point $\f$-dominates another, contradicting the definition of the Pareto front $\f X^*(\f)$.
\end{prf}

\begin{figure}[t]
\includegraphics[width=\hsize]{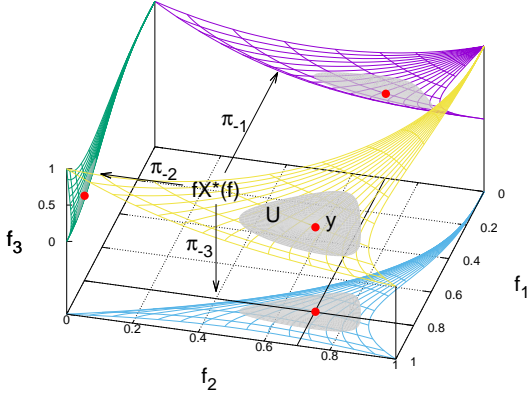}
\caption{A 3-objective (possibly non-simple) problem $\f = \set{f_1, f_2, f_3}$ and its Pareto front projections. The projection $\pi_{-i}$ restricted to $\Int \f X^*(\f)$ is injective and thus an embedding. By this fact, any interior point $y$ and its neighborhood $U$ in $\f X^*(\f)$ are mapped to an interior point $\pi_{-i}(y)$ and its neighborhood $\pi_{-i}(U)$ in $\pi_{-i} \f X^*(\f)$.}
\label{fig:projection}
\end{figure}
This lemma asserts that the interior of the Pareto front can be flattened while keeping its topology and ordering.
As shown in Fig.~\ref{fig:projection}, the projection induces a coordinate neighborhood $(\pi_{-i}, U)$ at any point $y \in \Int \f X^*(\f)$ where for any points $y, y' \in U$, it holds that $y_j \le y'_j \Leftrightarrow \pi_{-i}(y)_j \le \pi_{-i}(y')_j$ for all $j \ne i$.
This property is the key to investigate the interaction between the topology and the dominance-ordering on the interior of the Pareto front.
\begin{lem} \label{lem:neighbor}
Consider an $m$-objective (possibly non-simple) problem $\f$ whose Pareto front $\f X^*(\f)$ forms an $(m-1)$-manifold.
The following statement holds for any $y \in \Int \f X^*(\f)$:
any neighborhood $U \subseteq \f X^*(\f)$ of $y$ has a point that $\g$-dominates $y$ and is $(\f \setminus \g)$-dominated by $y$ for each $\g$ such that $\emptyset \subset \g \subset \f$.
\end{lem}
\begin{prf}
By Lemma~\ref{lem:projection}, the projection $\pi_{-i}$ restricted to $\f X^*(\f)$ is an embedding, and thus $\pi_{-i}(U)$ is a neighborhood of a point $\pi_{-i}(y)$ in $\R^{m-1}$.
Then, we can take an $(m-1)$-hyper-cube centered at $\pi_{-i}(y)$ in $\f X^*(\f)$ and can write its vertices as $(y_1 \pm \varepsilon, \dots, y_{i-1} \pm \varepsilon, y_{i+1} \pm \varepsilon, \dots, y_m \pm \varepsilon)$.
Here, $\varepsilon$ is a small positive number and signs, $\pm$, run over all possible combinations.
Among these vertices, let $\pi_{-i}(v)$ be the one such that coordinates related to $\g \subseteq \f_{-i} = \f \setminus \set{f_i}$ are $-\varepsilon$ and related to $\f_{-i} \setminus \g$ are $+\varepsilon$.
Thus, it holds that $y \dominated_{\g} v$ and $y \dominates_{(\f_{-i} \setminus \g)} v$.
Especially, $\g = \f_{-i}$ implies $y \dominates_{f_i} v$ and $\g = \emptyset$ implies $y \dominated_{f_i} v$.
Otherwise, it means that $y \dominates_{\f} v$ or $y \dominated_{\f} v$, contradicting $y, v \in \f X^*(\f)$.
Repeating the above argument for all $i$, we complete the proof.
\end{prf}

Other than the Pareto front $\f X^*(\f)$, do Lemmas~\ref{lem:projection} and~\ref{lem:neighbor} extend to the Pareto set image $\f X^*(\g)$ of a subproblem $\g \subproblem \f$?
The answer is NO: for general problems, the projection $\pi_{-i}$ is not an embedding of $\Int \f X^*(\g)$.
\begin{eg} \label{eg:non projection}
Again consider the problem in Example~\ref{eg:non pareto}:
\[
\minimize_{x \in [0,1]} f_1(x)=0, f_2(x)=x.
\]
The Pareto set of the subproblem $f_1$ is $X^*(f_1) = [0,1]$, and its image is $\f X^*(f_1) = \set{ (0, y) \mid 0 \le y \le 1}$.
The interior $\Int \f X^*(f_1) = \set{ (0, y) \mid 0 < y < 1}$ can be projected to $\pi_{-2} \Int \f X^*(f_1) = f_1 X^*(f_1) = \set{0}$, which implies $\pi_{-2}|_{\Int \f X^*(f_1)}$ is not injective and thus not an embedding.
Furthermore, since $X^*(f_1, f_2) = \set{0}$, we see that the set $X^*(f_1) \setminus X^*(f_1, f_2) = (0, 1]$ is weakly Pareto-optimal and not Pareto-optimal in $\f$.
The existence of these weak Pareto optima disrupts the injectivity of $\pi_{-2}$.
\end{eg}

Contrary, if the problem is simple, then there are no solutions that are weakly Pareto-optimal and not Pareto-optimal, which enable us to extend Lemmas~\ref{lem:projection} and~\ref{lem:neighbor} to the image $\f X^*(\g)$.
\begin{cor} \label{cor:projection simple}
Consider a simple problem $\f$ and a subproblem $\g \subproblemeq \f$.
The restriction of the projection $\pi_{-i}$ to $\Int \f X^*(\g)$ is an embedding.
\end{cor}
\begin{prf}
If $\pi_{-i}|_{\Int \f X^*(\g)}$ is not injective, then $\Int \f X^*(\g)$ contains two points that have the same coordinates except for the $i$-th value.
This means that $\f$ has a weak Pareto solution which is non-Pareto, contradicting Proposition~\ref{prp:no WP}.
\end{prf}

\begin{cor} \label{cor:neighbor}
Consider a simple problem $\f$ and a subproblem $\g \subproblemeq \f$.
For any $y \in \Int \f X^*(\g)$, the following statement holds:
for any neighborhood $U \subseteq \f X^*(\g)$ of $y$, there exists a point that $\h$-dominates $y$ and is $(\g \setminus \h)$-dominated by $y$ for each $\h$ such that $\emptyset \subset \h \subset \g$.
\end{cor}
\begin{prf}
Chose an objective function $f_i \in \f \setminus \g$ and let the remaining set be $\f_{-i} = \f \setminus \set{f_i}$.
By Corollary~\ref{cor:projection simple}, $\f X^*(\g)$ is mapped to $\f_{-i} X^*(\g)$ by $\pi_{-i}$, homeomorphically.
Next, chose another $f_j \in \f_{-i} \setminus \g$ and let the remainder be $\f_{-ij} = \f_{-i} \setminus \set{f_j}$.
Again by Corollary~\ref{cor:projection simple}, $\f_{-i} X^*(\g)$ is mapped to $\f_{-ij} X^*(\g)$ by $\pi_{-j}$, homeomorphically.
Though the repeated application of projections as long as it is an embedding, the original Pareto set image is finally mapped to $\g X^*(\g)$.
Let the composite of used projections $\pi_{-i}, \pi_{-j}, \dots$ be $\pi_{-*}$.
Generally, the composite of embeddings is again an embedding.
Thus, a point $y \in \Int \f X^*(\g)$ and its neighborhood $U \subseteq \f X^*(\g)$ is mapped homeomorphically to a point $\pi_{-*}(y) \in \Int \g X^*(\g)$ and its neighborhood $\pi_{-*}(U) \subseteq \g X^*(\g)$.
Proposition~\ref{prp:simple homeo} asserts $\g X^*(\g) \homeo \Delta^{\card{\g}-1}$, which together with Lemma~\ref{lem:neighbor} completes the proof.
\end{prf}

By this property, we can see the simplicity ensures that Pareto set images of subproblems are located on the boundary of their superproblems.
\begin{lem} \label{lem:boundary in F}
For a simple problem $\f$ and subproblems $\g,\h$ such that $\h \subset \g \subseteq \f$, the following relation holds:
\[
\boundary \f X^*(\g) \supseteq \f X^*(\h).
\]
\end{lem}
\begin{prf}
Suppose there exists a point $y \in \f X^*(\h)$ that is an interior point of $\f X^*(\g)$.
Now, $y \in \Int \f X^*(\g)$ and $\h \subset \g$ holds and thus by Corollary~\ref{cor:neighbor}, some neighborhood $U \subseteq \f X^*(\g)$ of $y$ has a point $\h$-dominating $y$.
This contradicts $y\in \f X^*(\h)$.
\end{prf}

Our next question is whether a similar relation holds for the Pareto set before mapped by $\f$.
To check this, we need the following lemma:
\begin{lem} \label{lem:comutative boundary}
For a simple problem $\f$ and a subproblem $\g \subproblemeq \f$, the map $\f$ commutes with the boundary and interior:
\begin{align}
\f \boundary X^*(\g) &= \boundary \f X^*(\g),\label{eqn:c1}\\
\f \Int X^*(\g) &= \Int \f X^*(\g).\label{eqn:c2}
\end{align}
\end{lem}
\begin{prf}
Generally, an embedding maps boundary to boundary and interior to interior.
Now $\f$ is an embedding on $X^*(\g) \subseteq X^*(\f)$, and thus it holds that $\f \boundary = \boundary \f$ and $\f \Int = \Int \f$.
\end{prf}

Using this fact, we show that under the simplicity, Pareto sets of subproblems have the same relation as its images.
\begin{cor} \label{cor:boundary in X}
For a simple problem $\f$ and a (proper!) subproblem $\g \subset \f$, the following relation holds:
\[
\boundary X^*(\f) \supseteq X^*(\g).
\]
\end{cor}
\begin{prf}
Since $\f$ is an embedding on $X^*(\f)$, there is the inverse map $\f^{-1}: \f X^*(\f) \to X^*(\f)$.
Lemma~\ref{lem:comutative boundary} converts Lemma~\ref{lem:boundary in F} to the assertion as follows:
\[
\begin{array}{crcl}
& \boundary \f X^*(\f) &\supseteq& \f X^*(\g)\\
\Leftrightarrow & \f \boundary X^*(\f) &\supseteq& \f X^*(\g)\\
\Leftrightarrow & \f^{-1} \f \boundary X^*(\f) &\supseteq& \f^{-1} \f X^*(\g)\\
\Leftrightarrow & \boundary X^*(\f) &\supseteq& X^*(\g).
\end{array}
\]
\end{prf}

The last key to the main theorem is the sphere embedding.
\begin{lem} \label{lem:sphere embedding}
Every embedding $f: S^n \to S^n$ is surjective and thus a homeomorphism where $S^n$ is an $n$-sphere $S^n = \set{x \in \R^{n+1} \mid \norm{x} = 1}$.
\end{lem}
\begin{prf}
Suppose $f$ is not surjective.
Then, there exists a point $y \in S^n \setminus f(S^n)$ and a stereographic projection with north pole $y$, denoted by $\pi: S^n \setminus \set{y} \to \R^n$.
Generally, any stereographic projection is an embedding, and the composite of embeddings is an embedding.
Therefore, $\pi f$ is an embedding of $S^n$ into $\R^n$.
This contradicts the well-known fact that $S^n$ cannot be embedded into $\R^n$.
\end{prf}
\begin{rem}
To keep the proof elementary, here we assumed that $S^n \not \embedding \R^n$ is known and derived that $S^n \embedding S^n$ is surjective.
There is an alternative proof deriving both by the Mayer-Vietoris exact sequence in a unified fashion.
Consult Hatcher~\cite{Hatcher02}, the two paragraphs after the proof of Proposition 2B.1 (pp.~169--170).
\end{rem}

Now, we show the goal of this section.
\begin{thm} \label{thm:assumption}
For a simple problem $\f$ and a subproblem $\g \subproblemeq \f$, it holds that
\begin{eqnarray}
\boundary X^*(\g)		&=& \bigsqcup_{\h \subset \g} \Int X^*(\h),\label{eqn:1}\\
\boundary \f X^*(\g)	&=& \bigsqcup_{\h \subset \g} \Int \f X^*(\h).\label{eqn:2}
\end{eqnarray}
\end{thm}
\begin{prf}
First, we show \eqref{eqn:2}.
What to be proven are:
\begin{description}
 \item[(a)] if $\h \ne \h'$, then $\Int \f X^*(\h) \cap \Int \f X^*(\h') = \emptyset,$
 \item[(b)] $\boundary \f X^*(\g) = \bigcup_{\h \subset \g} \Int \f X^*(\h)$.
\end{description}

(a)
When $\h \supset \h'$, it holds from Lemma~\ref{lem:boundary in F} that $\boundary \f X^*(\h) \supseteq \f X^*(\h')$.
The same holds for the inverse case $\h \subset \h'$.
We thus consider the case there is no inclusion relation between $\h$ and $\h'$.
Assume there exists a point $y \in \Int \f X^*(\h) \cap \Int \f X^*(\h')$, and let $U$ be a neighborhood of $y$ in $\f X^*(\h \cup \h')$.
If there exists a point $z \in U \setminus \paren{\f X^*(\h) \cup \f X^*(\h')}$, then $y \dominates_{\h} z$ and $y \dominates_{\h'} z$ hold, implying $y \dominates_{(\h \cup \h')} z$.
This contradicts $z \in U \subseteq \f X^*(\h \cup \h')$.
If $z$ does not exist, then the dimension of $U$ must be equal to that of $\f X^*(\h)$ or $\f X^*(\h')$, which contradicts $\dim U = \card{\h \cup \h'}-1$.
Consequently, such $y$ cannot exist.

(b)
Since Lemma~\ref{lem:boundary in F} states $\boundary \f X^*(\g) \supseteq \f X^*(\h)$, we have
\begin{align*}
\boundary \f X^*(\g) &\supseteq \bigcup_{\h \subset \g} \f X^*(\h)\\
&\supseteq \bigcup_{\h \subset \g} \Int \f X^*(\h).
\end{align*}
Thus, there is the inclusion map $\iota: \bigcup_{\h \subset \g} \Int \f X^*(\h) \embedding \boundary \f X^*(\g)$.
Generally, any inclusion map is an embedding.
Combining (a) with Lemma~\ref{lem:boundary in F}, we have $\f X^*(\h) \cap \f X^*(\h') = \f X^*(\h \cap \h')$.
Therefore, $\bigcup_{\h \subset \g} \Int \f X^*(\h)$ is the union of manifolds each homeomorphic to a simplex, which are glued as the faces of a simplex.
This fact ensures $\bigcup_{\h \subset \g} \Int \f X^*(\h) \homeo \boundary \Delta^{\card{\g}-1} \homeo S^{\card{\g}-2}$.
Contrary, $\f X^*(\g) \homeo \Delta^{\card{\g}-1}$ implies $\boundary \f X^*(\g) \homeo \boundary \Delta^{\card{\g}-1} \homeo S^{\card{\g}-2}$.
As Lemma~\ref{lem:sphere embedding} ensures that $S^n \embedding S^n$ is surjective, the inclusion map $\iota$ is surjective, implying that $\bigcup_{\h \subset \g} \Int \f X^*(\h) = \boundary \f X^*(\g)$ holds.

We can get \eqref{eqn:1} by converting \eqref{eqn:2} with Lamma~\ref{lem:comutative boundary}:
\[
\begin{array}{crcl}
& \boundary \f X^*(\g) &=& \bigsqcup_{\h \subset \g} \Int \f X^*(\h)\\
\Leftrightarrow & \f \boundary X^*(\g) &=& \bigsqcup_{\h \subset \g} \f \Int X^*(\h)\\
\Leftrightarrow & \f \boundary X^*(\g) &=& \f \bigsqcup_{\h \subset \g} \Int X^*(\h)\\
\Leftrightarrow & \f^{-1} \f \boundary X^*(\g) &=& \f^{-1} \f \bigsqcup_{\h \subset \g} \Int X^*(\h)\\
\Leftrightarrow & \boundary X^*(\g) &=& \bigsqcup_{\h \subset \g} \Int X^*(\h).
\end{array}
\]
From the second to third lines, we used the general property of a map $f(A) \cup f(B) = f(A \cup B)$.
\end{prf}


\subsection{Relation to Scalarization} \label{sec:AWA}
Equations \eqref{eqn:c1}--\eqref{eqn:2} together define a gluing structure of the Pareto sets and their images of subproblems of a simple problem.
This structure induces a natural stratification%
\footnote{The smoothness of the stratification is determined from that of the evaluation map.}
of the Pareto set (resp.\ the Pareto front) where each stratum is the interior of the Pareto set (resp.\ its image) of a subproblem.
Therefore, we can numerically compute the stratification by solving each subproblem.
Points spreading over all strata can be a good covering of the Pareto set/front.

To see why this structure enables decomposition-based EMO algorithms to cover the Pareto set/front, consider the weighted Tchebyshev-norm scalarization
\begin{equation} \label{eqn:scalarization}
\minimize_{x \in X} f_w(x) = \max_i w_i \paren{f_i(x) - z_i}
\end{equation}
where the weight $w=(w_1, \dots, w_m)$ is chosen from $\Delta^{m-1}$ and the utopian point is fixed to be $z_i = \inf_{x \in X} f_i(x)$.
Let $e_i$ be the $i$-th standard base in $\R^m$ whose $i$-th coordinate is one and the other coordinates are zero. 
The standard $(m-1)$-simplex is rewritten as $\Delta^{m-1} = [e_1, \dots, e_m]$.
Using the notation of the weight-optima correspondence
\[
S(W) = \bigcup_{w \in W} X^*(f_w),
\]
a well-known fact of the optima to \eqref{eqn:scalarization} can be written as
\begin{equation} \label{eqn:Chebyshev solution}
S([e_{i_1}, \dots, e_{i_k}]) = \wX(f_{i_1}, \dots, f_{i_k})
\end{equation}
for any choice of an arbitrary number of indices $i_1, \dots, i_k \in \set{1, \dots, m}$.
If the problem is simple, then we can go further: Proposition~\ref{prp:no WP} extends \eqref{eqn:Chebyshev solution} to
\[
S([e_{i_1}, \dots, e_{i_k}]) = X^*(f_{i_1}, \dots, f_{i_k}),
\]
and by Corollary~\ref{cor:boundary in X} we have
\[
S(\boundary [e_{i_1}, \dots, e_{i_k}]) = \boundary X^*(f_{i_1}, \dots, f_{i_k}).
\]
Therefore, a weight on each face gives a boundary point of each stratum with corresponding indices.

Unfortunately, the $L^\infty$-norm, as well as other existing scalarization methods including the weighted sum, the augmented Chebyshev-norm, PBI~\cite{Zhang07}, and IPBI~\cite{Sato14}, does NOT give the correspondence between the interiors:
\[
S(\Int [e_{i_1}, \dots, e_{i_k}]) \ne \Int X^*(f_{i_1}, \dots, f_{i_k}).
\]
Nevertheless, once boundary points of a stratum are obtained, we can find new weights corresponding to interior points of the stratum by interpolating the weights used for the boundary points.
Thus, the grid arrangement or divide-and-conquer generation of weights over $[e_{i_1}, \dots, e_{i_k}]$ practically often hit interior points of $X^*(f_{i_1}, \dots, f_{i_k})$.

\section{Simplicity of Benchmarks} \label{sec:benchmark}
This section investigates the simplicity of benchmark problems widely-used in the EMO community: ZDT suite~\cite{Zitzler00}, DTLZ suite~\cite{Deb05b}, WFG suite~\cite{Huband06}, and MED problem~\cite{Hamada11b}.

\subsection{ZDT Suite}
The ZDT suite has six $n$-variable 2-objective problems named ZDT1--6.
The decision variables are split into the \emph{position variables} $y = (y_1, \dots, y_k)$ $(0 < k < n)$ and the \emph{distance variables} $z = (z_1, \dots, z_l)$ $(l = n - k)$, defining the problems in the following unified format:
\begin{equation*}
\begin{split}
\minimize_{(y, z) \in Y \times Z}  \f(y,z) &= (f_1(y,z), f_2(y,z)),\\
\text{where }
f_1(y,z) &= f(y_1),\\
f_2(y,z) &= g(z) h(f_1(y), g(z)),\\
X = Y \times Z &=(Y_1 \times \dots \times Y_k) \times (Z_1 \times \dots \times Z_l),\\
Y_1 &= [0,1],\\
Y_2, \dots, Y_k, Z_1, \dots, Z_l &=
\begin{cases}
 [-5,5] & \text{(ZDT4)},\\
 [0,1] & \text{(otherwise)}.
\end{cases}
\end{split}
\end{equation*}
Users can make different problems by changing placeholder functions $f, g, h$.
For the concrete specification of $f, g, h$ for ZDT1--6, see Zitzler et~al.~\cite{Zitzler00}.
The above general formulas are enough to show that the problems are non-simple.

\begin{thm} \label{thm:ZDT}
ZDT1--6 are all non-simple, independent of the choice of variable dimension $n$ and position-variable dimension $k$.
Additionally, this suite cannot create simple problems no matter how $f,g,h$ are specified unless their domains are modified.
\end{thm}
\begin{prf}
First, we exclude ZDT5 from the following analysis since it is a binary-variable problem which is clearly non-simple.
Then, for ZDT1--4, 6, the function $f$ defining $f_1$ depends on a single variable, $y_1$.
The other variables can take an arbitrary value on the optima of $f$, and thus $X^*(f_1) = Y_1^*(f) \times Y_2 \times \dots \times Y_k \times Z$.
This means $X^*(f_1) \not \homeo \Delta^0$ and contradicts the simplicity condition (S1).
Consequently, ZDT1--6 are all non-simple.
\end{prf}

The reason why this suite cannot be simple is that $f$ depends only on $y_1$.
Generally, when the problem has an objective function independent of some variables, its Pareto set extends to higher dimensions than usual, contradicting the simplicity condition (S1).
The existence of unused variables is a quick test for non-simplicity.

\subsection{DTLZ Suite}
The DTLZ suite consists of nine problems named DTLZ1--9.
Their decision variables are split into position variables $y$ and distance variables $z$, as ZDT, but the number of objectives $m$ can be set arbitrarily.
See Deb et~al.~\cite{Deb05b} for definition.

\begin{thm} \label{thm:DTLZ}
DTLZ1--9 are all non-simple, independent of the choice of variable dimension $n$, objective dimension $m$, and position-variable dimension $k$.
\end{thm}
\begin{prf}
Every problem has an objective function ignoring some variables as follows:
\begin{eqnarray*}
\text{DTLZ1--6: } & f_m(y,z)&=(1+g(z))f(y_1)\\
\text{DTLZ7: }      & f_1(y,z)&=y_1\\
\text{DTLZ8, 9: }  & f_1(y,z)&=y_1^{0.1} + y_2^{0.1} + \dots + y_{\lfloor n/m \rfloor}^{0.1}
\end{eqnarray*}
Therefore, DTLZ1--9 are all non-simple.
\end{prf}

Furthermore, Huband et~al.~\cite{Huband06} Table VII shows that DTLZ1--6 has $\f$ that is not injective on $X^*(\f)$ and DTLZ7 has a disconnected Pareto front.
This is another evidence for the non-simplicity of DTLZ1--7.

\subsection{WFG Suite}
The WFG suite contains nine problems, WFG1--9, having the form:
\begin{equation*} \label{eqn:WFG}
\begin{split}
\minimize_{(y,z) \in Y \times Z}  \f(x) &= (f_1(x), \dots, f_m(x)),\\
\text{where } f_i(x) &= x_m + 2i \times h_i(x_1, \dots, x_{m-1}),\\
 t_Y &:\left \{
  \begin{array}{ccc}
   Y &\to& X_1 \times \dots \times X_{m-1},\\
   (y_1, \dots, y_k) &\mapsto& (x_1,\dots,x_{m-1}),
  \end{array}\right.\\
 t_Z &:\left \{
  \begin{array}{ccc}
   Z &\to& X_m,\\
   (z_1, \dots, z_l) &\mapsto& x_m,
  \end{array}\right.\\
 X &= X_1 \times \dots \times X_m = [0,1]^m,\\
 Y &= [0,2] \times [0,4] \times \dots \times [0,2k],\\
 Z &= [0,2k+2] \times [0,2k+4] \times \dots \times [0,2n].
\end{split}
\end{equation*}
The functions $h_1, \dots, h_m, t_Y, t_Z$ are placeholders.
In this suite, the variables $y$ and $z$ are mapped by the \emph{transformation functions} $t_Y$ and $t_Z$ to the \emph{position variables} $x_1, \dots, x_{m-1}$ and the \emph{distance variable} $x_m$, then passed to the objective functions $f_i$.
For this reason, $y$ and $z$ are called the \emph{position-related variables} and the \emph{distance-related variables}, respectively.
For the concrete definition, see Huband et~al.~\cite{Huband06}.

\begin{thm} \label{thm:WFG}
WFG2, 4, 5, 9 are always non-simple.
WFG1, 3, 6--8 are simple if and only if the dimension of the position-related variables $y$ is $k=1$.
Here, one can set $k=1$ only when the number of objectives is $m=2$ because Huband et~al.~\cite{Huband06} Table~XIV shows that these problems require $k \bmod (m-1) = 0$.
\end{thm}
\begin{prf}
First, consider WFG2.
Huband et~al.~\cite{Huband06} Table~XIV shows that WFG2 has a disconnected Pareto front.
Such a front cannot be homeomorphic to $\Delta^{m-1}$, which contradicts the property of a simple problem shown in Proposition~\ref{prp:simple homeo}.
The following discussion treats the rest of the problems.

Let us check the properties of the Pareto set $X^*(\f)$ in the transformed variable space $X = [0, 1]^m$ and the map $\f$ on $X^*(\f)$.
As described in Huband et~al.~\cite{Huband06}, all problems%
\footnote{As opposed to Huband et~al.~\cite{Huband06}, WFG2 actually has a different Pareto set in which $x_1$ is conditioned to be Pareto-optimal.
To avoid a complication caused by this difference, we first finished WFG2.}
have
\[
X^*(\f) = [0,1]^{m-1} \times \set{0}.
\]
By the properties of the \emph{shape functions} $h_i$ shown in Huband et~al.~\cite{Huband06} Table~X, $\f$ is an embedding upon $X^*(\f)$ if and only if $m=2$.%
\footnote{Except for \textbf{Disconnected} used only in WFG2.}
For $m > 2$, it holds that
\[
\f(0, x_2, \dots, x_{m-1}, 0) = 0,
\]
and thus $\f$ is not injective on $X^*(\f)$.

Since the transformation functions, $t_Y$ and $t_Z$, are surjective as described in Huband et~al.~\cite{Huband06}, the composite evaluation map $\f t: Y \times Z \xrightarrow{(t_Y, t_Z)} X \xrightarrow{\f} \R^m$ is not injective on the Pareto set $(Y \times Z)^*(\f t)$ in the untransformed variable space $Y \times Z$ for $m > 2$.
This does not meet the simplicity condition (S2).

In the case of $k > m-1$, the transformation function $t_Y$ involves \textbf{Reduction: Weighted Sum} or \textbf{Reduction: Non-separable} to decrease the dimension of $Y$.
By examining Huband et~al.~\cite{Huband06} Table~XI, we can see that both functions are not injective.
Therefore, $t_Y$ maps two different points $y, y' \in Y$ to the same Pareto solution $(x_1,\dots, x_{m-1}) \in [0,1]^{m-1}$, contradicting  the simplicity condition (S2).

There remains the case $k=1$ and $m=2$ where the problems may be simple.
First, let us consider the 1-objective subproblems.
From Huband et~al.~\cite{Huband06} Table~X, we have
\[
X^*(f_1)=\set{(0,0)},\quad X^*(f_2)=\set{(1,0)}.
\]
This implies that the simplicity condition (S1) for these problems is equivalent to the following criterion:
\[
t_Y^{-1}(0), t_Y^{-1}(1), t_Z^{-1}(0) \text{ are all a point}.
\]
Since (S2) automatically follows from (S1) when the problem is 1-objective, the above criterion is a necessary and sufficient condition for the simplicity.
We can see from Huband et~al.~\cite{Huband06} Table~XI that $t_Y^{-1}(0)$ and $t_Z^{-1}(0)$ are always a point, but it depends on the case whether $t_Y^{-1}(1)$ is a point or not.
WFG4, 5, 9 introduce \textbf{Shift: Deceptive} or \textbf{Shift: Multi-modal} into $t_Y$, making $t_Y^{-1}(1)$ not a point.
Thus, these problems are non-simple.
For WFG1, 3, 6--8, $t_Y^{-1}(1)$ becomes a point, and their 1-objective subproblems are simple.
Next, let us consider their 2-objective subproblems.
It holds that
\[
X^*(f_1, f_2) = [0,1] \times \set{0},
\]
and it has been confirmed that $\f = \set{f_1, f_2}$ is an embedding on $X^*(f_1, f_2)$ and $t_Z^{-1}(0)$ is a point.
Therefore, the equivalent condition to the simplicity is as follows:
\[
t_Y^{-1}: [0,1] \to [0, 2] \text{ embedding}.
\]
For $k=1$ and $m=2$, the transformation function $t_Y$ in WFG1, 3, 6--8 can be simplified into the form $t_Y(y)=y^\alpha$.
Thus, its inverse, $t_Y^{-1}(x)=x^{1/\alpha}$, is an embedding.
Now, we have checked that WFG1, 3, 6--8 are simple if and only if $k=1$ and $m=2$.
\end{prf}

Note that Huband et~al.~\cite{Huband06} Table~XIV describes that WFG3 has a degenerate Pareto front, which seems to be an evidence that WFG3 is always non-simple.
However, the degeneracy actually occurs only when $m > 2$.
In our analysis for $m=2$, the Pareto front of WFG3 forms a line segment, which does not disrupt the simplicity.

\subsection{MED}
This is a single problem, MED, defined as follows:
\begin{equation} \label{eqn:MED}
\begin{split}
\minimize_{x \in X=\R^{n}}  \f(x) &= (f_1(x), \dots, f_m(x)),\\
\text{where }	 f_i(x) &= \norm{x - x_i^*}^{p_i},\\
		 x_i^* &= (\underbrace{0, \dots, 0,}_{i-1} 1 \underbrace{, 0, \dots, 0}_{n-i}),\\
		 0 & < p_i < \infty.
\end{split}
\end{equation}
The \emph{front-shape parameters}, $p_i$, which determine the convexity of the Pareto front, are user-specified parameters as well as the variable dimension $n$ and the number of objectives $m$.

\begin{thm} \label{thm:MED}
MED is always simple independent of the choice of parameters $n, m, p_i$.
Additionally, changing individual optima, $x_i^*$, does not break the simplicity as long as they are affinely independent.
\end{thm}
\begin{prf}
First, consider the case of $p_i = 1$.
This corresponds to a facility location problem under the $L^2$-norm.
The Pareto set of this problem is known as the convex hull of $x^*_1, \dots, x^*_m$~\cite{Ward89}.
Thus, if $m \le n+1$ holds and $x^*_1, \dots, x^*_m$ are affinely independent, then the convex hull is the $(m-1)$-simplex spanned by $x^*_1, \dots, x^*_m$.
Indeed, by the definition of $x^*_i$, this problem can be defined only when $m \le n$, and $x^*_1, \dots, x^*_m$ are affinely independent.
Thus, the Pareto set $X^*(\f)$ is an $(m-1)$-simplex, which ensures that the problem $\f$ satisfies the simplicity condition (S1).
Analyzing the gradient of the $L^2$-norm, we can see that the map $\f$ is an embedding on $X^*(\f)$, which satisfies (S2).
The same argument applies to the subproblems, confirming that they satisfy (S1) and (S2); the problem $\f$ is simple.

The case $p_i \ne 1$ can be considered as the composite of $f_i$ in the case of $p_i = 1$ and the $p_i$-th power.
Since any positive power is an order-preserving homeomorphism $[0, \infty) \to [0, \infty)$, the composition preserves the simplicity of the facility location problem.
\end{prf}

\section{Conclusions} \label{sec:conclusions}
In this paper, we have discussed the simple problem and showed that the Pareto sets of its subproblems (resp.\ their images) constitute a stratification of its Pareto set (resp.\ its Pareto front).
This topological property gives a theoretical guarantee that decomposition-based EMO algorithms can obtain an entire approximation of the Pareto set as well as the Pareto front.
We have also investigated the simplicity of benchmark problems widely-used in the EMO community.
All problems in the ZDT and DTLZ suites are non-simple.
The WFG suite contains five simple problems under a very restrictive situation but usually does not, whereas the MED problem is always simple.

We believe that the absence of simple problems in the standard benchmark suites is a considerable gap between the benchmark and the real-world since there are many evidences that a large portion of nowadays applications seems to be simple.
Additionally, real-world applications involving simulations can be black-box; it would be important to develop an estimation method for the simplicity of black-box problems from a finite set of approximate solutions.

\bibliographystyle{ACM-Reference-Format}
\bibliography{ref} 

\end{document}